\documentclass[reqno,12pt,twoside,english]{amsart}
\usepackage[utf8]{inputenc}
\usepackage[T1]{fontenc}
\usepackage{amsmath, amsthm}
\usepackage{amsfonts, amssymb}
\usepackage{enumerate}
\usepackage{graphicx}
\usepackage{xcolor}
\usepackage{float}
\usepackage[a4paper,bindingoffset=0.13in,%
left=1in,right=1in,top=0.6in,bottom=0.6in,%
footskip=.15in]{geometry}
\usepackage{comment}
\usepackage{enumitem}
\usepackage{amscd}
\usepackage{subfiles}
\usepackage{tikz}

\newcommand*{\R}{\mathbb{R}}
\newcommand*{\closure}[1]{\overline{#1}}
\newcommand*{\bdary}[1]{\partial #1}
\newcommand*{\dR}{\ensuremath{\dot{\mathbb{R}}^n}}
\newcommand*{\Sf}{\mathbb{S}}
\newcommand*{\B}{\mathbb{B}}
\newcommand*{\D}{\mathbb{D}}
\newcommand*{\eps}{\varepsilon}

\DeclareMathOperator{\diam}{diam}

\numberwithin{equation}{section}
\theoremstyle{definition}		
\newtheorem{theorem}{Theorem}[section]
\newtheorem{proposition}[theorem]{Proposition}
\newtheorem{lemma}[theorem]{Lemma}

\newtheorem{example}[theorem]{Example}

\title[Characterizations of generalized John domains]{Characterizations of generalized John domains via homological bounded turning}
\author[P. Goldstein]{Pawe\l{} Goldstein}
\address[Pawe\l{} Goldstein]{Institute of Mathematics, Faculty of Mathematics, Informatics and Mechanics, University of Warsaw,
	Banacha 2, 02-097 Warsaw, Poland}
\email{goldie@mimuw.edu.pl}

\author[Z. Grochulska]{Zofia Grochulska}
\address[Zofia Grochulska]{Faculty of Mathematics, Informatics and Mechanics, University of Warsaw,
	Banacha 2, 02-097 Warsaw, Poland}
\email{z.grochulska@uw.edu.pl}

\author[C.-Y. Guo]{Chang-Yu Guo}
\address[Chang-Yu Guo]{Research Center for Mathematics and Interdisciplinary Sciences, Shandong University, 266237 Qingdao, P. R. China and  Frontiers Science Center for Nonlinear Expectations, Ministry of Education, P. R. China}
\email{changyu.guo@sdu.edu.cn}

\author[P. Koskela]{Pekka Koskela}
\address[Pekka Koskela]{Department of Mathematics and Statistics, University of Jyv\"askyl\"a, P.O. Box 35, FI-40014 Jyv\"askyl\"a, Finland}
\email{pkoskela@maths.jyu.fi}

\author[D. Nandi]{Debanjan Nandi}
\address[Debanjan Nandi]{Department of Mathematics, the Weizmann Institute, Israel.}
\email{debanjan.nandi@weizmann.ac.il}

\thanks{P. Goldstein was partially supported by FNP grant POMOST BIS/2012-6/3 and by NCN grant no 2012/05/E/ST1/03232. C.-Y. Guo is supported by the Young Scientist Program of the Ministry of Science and Technology of China (No.~2021YFA1002200), the National Natural Science Foundation of China (No.~12101362) and the Natural Science Foundation of Shandong Province (No.~ZR2021QA003). P. Koskela was partially supported by the Academy of Finland grant 323960.}

\begin{document}
	
	\begin{abstract}
		In this paper, we extend the characterization of John disks obtained by N\"akki and V\"ais\"al\"a \cite{nv91} to generalized John domains in higher dimensions under mild assumptions.
		%Simple examples indicate that our assumptions for such a characterization are optimal.
		The main ingredient in this characterization is to use the higher dimensional analogues of the local linear connectivity (LLC) and homological bounded turning properties introduced by V\"ais\"al\"a in his study of metric duality theory \cite{v97}.
		
		Somewhat surprisingly, we constructed a uniform domain in $\R^3$, which is topologically simple, such that the complementary domain fails to be homotopically $1$-bounded turning. In particular, this shows that a similar characterization of generalized John domains in terms of higher dimensional homotopic bounded turning does not hold in dimension three.
	\end{abstract}
	
	\maketitle
	
\section{Introduction and main results}

One of the most important and ubiquitous notions in the theory of domains and its applications to modern analysis is that of a \emph{John domain.}
Recall that an open, connected and bounded set $\Omega\subset \R^n$ is a \emph{John domain} if there is  a point $x_0\in \Omega$ and a~constant $C$ such that each $x\in \Omega$ can be connected with $x_0$ with a rectifiable \emph{John curve} $\gamma: [0,1]\to \Omega$, $\gamma(0)=x$,
$\gamma(1)=x_0$, satisfying for all $t \in [0,1]$
\begin{equation}\label{eq:def for John}
	l(\gamma([0,t])) \leq Cd(\gamma(t), \partial \Omega).
\end{equation}

The term \emph{John domain} was introduced by Martio and Sarvas~\cite{ms79} to honour F. John, who used this condition in his seminal work on elasticity~\cite{j61}. The class of John domains includes all smooth domains, Lipschitz domains and certain fractal domains, such as the interior of the von Koch snowflake and other snowflake-type domains.

In the planar case, N\"akki and V\"ais\"al\"a~\cite{nv91} developed a rich theory of \emph{John disks}, i.e., simply connected planar John domains and proved the following important characterization of John disks (for definitions see Section \ref{sec: definitions}). The exact formulation contains a~few further equivalencies, which are however not in the scope of this paper.

\begin{proposition}\cite[Theorem 4.5]{nv91} \label{thm:Nakki and Vaisala} 
	Let $\Omega\subset \dot{\mathbb{R}}^2$ be a simply connected domain such that $\overline{\Omega}$ does not contain the point at infinity. Then the following conditions are quantitatively equivalent:
	\begin{itemize}
		\item[1)] $\Omega$ is a John disk,
		\item[2)] $\Omega$ is LLC-2,
		\item[3)] $\dot{\mathbb{R}}^2 \setminus \closure{\Omega}$ is bounded turning.
	\end{itemize}
\end{proposition}

In \cite[Theorem 4.5]{nv91}, the theorem is stated for conformal disks. By 
applying the Riemann mapping theorem, it thus holds for any simply connected domain $\Omega\subset \dot{\mathbb{R}}^2$ such that  $\overline{\Omega}$ does not contain the point at infinity.

The theory of John domains has found many applications in geometry and analysis, in particular in the area of geometric function theory, starting from the beautiful work of Ahlfors \cite{a63}, who proved that a bounded domain $\Omega\subset \R^2$ is a quasidisk (image of the unit disk under a global planar quasiconformal mapping) if and only if both $\Omega$ and its complementary domain are LLC-2. This is further equivalent via Theorem \ref{thm:Nakki and Vaisala} to the statement that both $\Omega$ and its complement are John disks. Naturally, the theory of quasidisks plays an important role in the theory of planar quasiconformal mappings, see \cite{AIM09,GH12}.

Recently, there has been a~considerable interest in developing the theory of mappings of finite distortion with locally exponentially integrable distortion (see \cite{IM01,HK14}), which is a~substantially larger class than quasiconformal mappings. Accordingly, one can define \emph{generalized quasidisks} to be images of the unit disk under such mappings. The results of \cite{gkt14, g15b} show that a generalized quasidisk is not necessarily John, but that it satisfies a nonlinear version of the John condition \eqref{eq:def for John}. To capture this nonlinear John property, Guo and Koskela introduced in~\cite{gk14} the class of \emph{$\varphi$-John domains} for any continuous, increasing $\varphi\colon [0,\infty)\to [0,\infty)$ with $\varphi(0)=0$ and $\varphi (r)\ge r$ for all $r>0$. This clearly is a~broader class than the class of John domains. Nonetheless, the quantitative equivalence of 1) and 2) in Theorem \ref{thm:Nakki and Vaisala} does admit a natural extension in the category of generalized John disks, see~\cite{gk14}.

In this paper, we present a generalization of Theorem \ref{thm:Nakki and Vaisala} to higher dimensions in the context of generalized John domains. Observe that now $\Omega$ is not assumed to be bounded since the notion of $\varphi$-John domains can be adapted to work for the unbounded case as well, see Section \ref{sec: definitions}.
\begin{theorem}\label{thm main: characterization}
	Let $\Omega \subsetneq \R^n$ with $H^1(\closure{\Omega})=0$ be a proper subdomain which satisfies the ball separation property and is locally collared along $\partial\Omega$. Then the following conditions are quantitatively equivalent:
	\begin{itemize}
		\item[1)] $\Omega$ is (0,$\varphi$)-John;
		\item[2)] $\Omega$ is $\varphi$-LC-2;
		\item[3)] the complementary domain $U=\dR\backslash \closure{\Omega}$ is homologically $(n-2,\varphi)$-bounded turning.
	\end{itemize}
\end{theorem}
The essence of this theorem is addressing a~fundamental question about the appropriate higher dimensional analogues of a~simply connected planar domain and the bounded turning property. In the planar case, due to the uniformization theorem, the fact that a (planar) domain is simply connected carries not only topological, but also geometric information. In higher dimensions simple examples show that an LLC-2 homeomorphic image of the unit ball in $\R^n$, $n\geq 3$, can fail to be John (see for instance Example~\ref{example:easy one}). Indeed, the proof of Theorem \ref{thm:Nakki and Vaisala} in  \cite{nv91} relies in a crucial way on the theory of planar conformal mappings.  Therefore, a trivial topological extension is not sufficient any more in higher dimensions and one expects a suitable geometric extension of planar simple connectedness. 

We propose as such the \emph{ball separation property} introduced in \cite{BuckleyKoskela} and studied further for instance in \cite{bhk01,BaloghBuckley}. We say that a~proper domain $\Omega\subset\dR$ has the \emph{ball separation property} if there exists a positive constant $C$ such that for each pair of points $x,y\in \Omega$, for any quasihyperbolic geodesic $[x,y]$ from $x$ to $y$ each $z\in [x,y]$, and for every curve $\gamma$ in $\Omega$ joining $x$ with $y$ we have
\begin{equation}\label{eq:separation property}
	B(z,Cd(z,\bdary\Omega))\cap \gamma\neq \emptyset.
\end{equation} 
This property has turned out to be useful in many problems, for example in connection with Sobolev-Poincar\'e inequalities \cite{BuckleyKoskela} and uniform continuity of quasiconformal mappings into irregular domains~\cite{g15}. It is also relatively common -- every Gromov hyperbolic domain satisfies the ball separation property \cite{BaloghBuckley}, which in particular implies that uniform domains (being Gromov hyperbolic) satisfy it as well \cite{bhk01}.

% For $p=0$ it is the classical definition (when was it introduced); What does $|f|$ mean?

As for the higher dimensional analogues of bounded turning property, we shall adapt for the nonlinear case the concepts of $p$-dimensional homological and homotopical bounded turning properties introduced by Alestalo and V\"ais\"al\"a in \cite{av94} and \cite{av96} (for definitions see Section \ref{sec: definitions}). Both definitions agree for $p = 0$ and are quantitatively equivalent to the classical bounded turning property, which states that any two points $x, y \in \Omega$ can be joined by a~continuum whose diameter does not exceed $C|x-y|$. The idea of $p$-dimensional bounded turning properties is closely related with the concept of \emph{homological joinability} introduced almost simultaneously by V\"ais\"al\"a in his beautiful work \cite{v97}. In fact, the key step in the proof of Theorem \ref{thm main: characterization} is the use of his Metric Duality Theorem (Theorem 2.7 in \cite{v97}), whose nonlinear version is stated in the Section \ref{sec: definitions}.

In general, homological bounded turning is a~weaker property than its homotopical counterpart; see \cite{a94,av96,av97} for more on the relation between homologically defined objects and homotopically defined objects. In particular, it was an open problem (see \cite[Open problem 4.8]{av94}) to determine whether for a~$p$-cell in $\R^n$ these two concepts are quantitatively equivalent.  As a partial answer to this problem, we show in Section \ref{sec: main example}, that in the case when $\Omega$ is a~planar, simply connected domain, these two concepts are indeed quantitatively equivalent for $p=1$. Thus one could still expect that a stronger form of Theorem~\ref{thm main: characterization} in terms of homotopic bounded turning could hold, if $\Omega$ enjoys a~suitable higher dimensional analogue of simple connectedness. Somewhat surprisingly, our main result  of this paper shows that this expectation is false when $n=3$, even in the linear case.
\begin{theorem}\label{thm main: example}
	There exists a domain $\Omega\subset \R^3$ with the following properties
	\begin{enumerate}
		\item $\closure{\Omega}$ is homeomorphic to the closed unit ball $\closure{\mathbb{B}}$;
		\item $\Omega$ is a uniform domain; 
		\item $U=\R^3\backslash \closure{\Omega}$ is homologically $(1,C)$-bounded turning for some constant $C>1$;
		\item $U$ is not homotopically $(1,C)$-bounded turning for any $C>1$.
	\end{enumerate}
\end{theorem}
Due to the fact that homeomorphisms between Euclidean spaces of the same dimension map boundaries onto boundaries and interiors onto interiors, $\Omega$ has to be homeomorphic to the open unit ball. Moreover, property (1) guarantees that $H_1(\closure{\Omega}) = 0$ and that $\Omega$ is locally collared along $\partial \Omega$ whereas property (2) implies that $\Omega$ satisfies the ball separation property. Thus it satisfies the assumptions of Theorem \ref{thm main: characterization}, in light of which we can conclude that property (3) follows from the fact that uniform domains are LLC-2.

Interestingly, the domain whose existence is claimed in Theorem~\ref{thm main: example} is the interior of an \emph{Alexander's horned ball} constructed so as to warrant its nice metric properties. The construction relies heavily on three dimensional topology, in particular, the close relation of the fundamental group and the first homology group. Thus, this kind of constructions do not generalize easily to higher dimensions (i.e., $n\geq 4$), when we deal with higher-dimensional homology groups. Then, the relation between homology and homotopy groups in the same dimension is far more complex. Therefore, even though we show in Theorem~\ref{thm main: example}, by the constructed counterexample, that a variant of Theorem~\ref{thm main: characterization} with homological bounded turning replaced by homotopic bounded turning assumption fails in dimension 3, we do not know whether it holds for $n\geq 4$.

\subsection*{Notation and structure}
The closure of a set $U\subset\R^n$ is denoted $\closure{U}$ and the boundary $\bdary{U}$. The open ball of radius $r>0$ centered at $x\in\R^n$ is denoted by $B(x,r)$ and in the case of the unit ball we write $\mathbb{B}$. The boundary of $B(x,r)$ will be denoted by $S(x,r)$ and in case of the unit sphere by $\mathbb{S}$. The symbol~$\Omega$ always refers to a domain, i.e., a~connected and open subset of~$\R^n$ or $\dR$. Whenever we write $\gamma(x,y)$ or $\gamma_{xy}$, it refers to a curve or an arc from $x$ to $y$. Given a pair of points $x,y\in \Omega$, we denote by $[x,y]$ a~quasihyperbolic geodesic that joins $x$ and $y$.

We denote by $\dR$ the one-point compactification of $\R^n$, that is $\dR=\R^n\cup\{\infty\}$ with distance $|\cdot|$ using the convention that $|\infty - \infty| = 0$ and $|a - \infty| = \infty$ for any $a \in \R^n$. The topology in $\dR$ is given by the stereographic projection homeomorphism with $\mathbb{S}^{n-1}$. For an open or closed set $X$ in $\dR$, we denote by $H_p(X)$ the reduced singular $p$-homology group of $X$ and by $H^p(X)$ the Alexander-Spanier $p$-cohomology group, both with coefficients in $\mathbb{Z}$. We say that condition $P$ implies condition $P'$ quantitatively if $P$ with constant $c$ implies $P'$ with some other constant $c'(c)$. If these conditions involve a~function $\varphi$, we say that $P$ implies $P'$ quantitatively if $P$ with function $\varphi$ implies $P'$ with function $\psi(r) = a \varphi(br) + cr$ for some positive constants $a, b, c$. We say that $P$ and $P'$ are quantitatively equivalent if $P$ implies $P'$ quantitatively and if $P'$ implies $P$ quantitatively as well.

The paper is organized as follows. Section~\ref{sec: definitions} fixes the notation and basic definitions. In Section~\ref{sec: main characterization}, we prove Theorem~\ref{thm main: characterization} and show its sharpness in terms of given assumptions. We also give some remarks on quasihyperbolic geodesics in domains with the ball separation property. Section \ref{sec: main example} is devoted to constructing an example satisfying Theorem \ref{thm main: example}. We also include an appendix to collect the necessary auxiliary results on generalized joinability.

\section{Preliminaries and Definitions} \label{sec: definitions}

Let $\varphi: [0, \infty) \to [0, \infty)$ be a~continuous, increasing function such that $\varphi(0) = 0$ and $\varphi(r) \geq r$. Whenever we refer to a~nonlinear definition we do so with respect to such a~function.

\subsubsection*{John domain}
We say that $\Omega \subset \dR$ is a~\emph{$(p,\varphi)$-John} domain if $\Omega$ is a~domain and if for every $p$-cycle $z \subset \Omega$ there is a $(p+1)$-chain $g\subset\Omega$ such that $\partial g=z$ and
\begin{equation} \label{eq: general John}
	d(x,|z|)\leq \varphi(d(x,\partial \Omega)) \quad \text{for all}\ x\in |g|.
\end{equation}
In particular, $H_p(\Omega) = 0$.

For $p=0$ and $\varphi(r) = Cr$ it means that any two points $x_1, x_2$ in $\Omega \subset \dR$ can be joined with a curve $\gamma$ connecting them so that
\begin{equation*}
	\min_{i=1,2} \{ |\gamma(t) - x_i| \} \leq C d(x, \Omega) \text{ for all } t.
\end{equation*}
This is the so called \emph{distance cigar} condition used for example in \cite{Vaisala_uniform} and \cite{nv91}. If $\infty \notin \overline{\Omega}$ (i.e., if $\Omega$ is a~bounded domain in $\R^n$), distance cigar condition is quantitatively equivalent to the original condition \eqref{eq:def for John} quoted in the Introduction (sometimes called the \emph{length-John} or \emph{length carrot} condition). It also quantitatively coincides with the notion of a~\emph{diam-John} domain with $\diam{\gamma([0,t])}$ instead of $l(\gamma([0,t]))$ on the left hand-side of \eqref{eq:def for John}, which was firstly proved in \cite[Theorem 2.7]{ms79}. Section 2 of \cite{nv91} provides a~nice overview of the interplay between different cigar and carrots conditions.

As in \cite{gk14}, we say that a bounded domain $\Omega$ in $\R^n$ is a \emph{$\varphi$-length John domain} if it satisfies the following counterpart of the length John condition \eqref{eq:def for John}:
\begin{equation*}
	l(\gamma([0,t])) \leq \varphi(Cd(\gamma(t),\bdary\Omega)).
\end{equation*}
Replacing $l (\gamma([0,t]))$ with $\diam \gamma([0,t])$ or $|\gamma(t) - x_0|$ results in the concepts of $\varphi$-dist and $\varphi$-diam John domains, respectively. For bounded domains, $p=0$ and general $\varphi$ definition \eqref{eq: general John} agrees with $\varphi$-dist John domains. It has been proved in \cite{gk14} that under some technical assumptions for $\varphi$ for planar simple connected domains being a $\varphi$-dist John domain is equivalent with being a $\varphi$-diam John domain. However, even in this setting, being a $\varphi$-length John domain for general $\varphi$ is not necessarily quantitatively equivalent with being a $\varphi$-diam John domain (although it is so in the linear case).

\subsubsection*{Local connectivity}
We say that a set $E\subset \dR$ is \emph{linearly locally connected} (LLC) if there exists $C>1$ such that for any $x \in \R^n$ and $r>0$
\begin{itemize}[leftmargin=2cm]
	\item[(LLC-1)] each pair of points in $B(x,r)\cap E$ can be joined in $B(x,Cr)\cap E$,
	\item[(LLC-2)] each pair of points in $E\setminus B(x,Cr)$ can be joined in $E\setminus B(x,r)$.
\end{itemize}
These conditions can be generalized to the nonlinear case as follows: $E$ is $(\varphi, \psi)$-LC if there exists $C>1$ such that for any $x \in \R^n$ and $r>0$
\begin{itemize}[leftmargin=2cm]
	\item[($\varphi$-LC-1)] each pair of points in $B(x,r)\cap E$ can be joined in
	$B(x,\varphi(r))\cap E$,
	\item[($\psi$-LC-2)] each pair of points in $E\backslash B(x, \psi(r))$ can be joined in $E\setminus B(x, r)$.
\end{itemize}
Depending on what is meant by joining, one can consider \emph{pathwise} and \emph{continuumwise} versions of these conditions. However, they coincide in case of domains in $\dR$ since if $E$ is locally compact and locally path-connected, then pathwise connectivity is quantitatively equivalent to continuumwise connectivity (see e.g.~\cite{hy88}).

\subsubsection*{Quasihyperbolic metric and quasihyperbolic geodesics}
The quasihyperbolic length of a~curve $\gamma$ joining point $x$ and $y$ in $\Omega \subsetneq \R^n$ is defined as
\begin{equation*}
	k_\Omega\text{-length}(\gamma)=\int_\gamma \frac{ds}{d(x,\bdary\Omega)}.
\end{equation*}
Then, the quasihyperbolic metric $k_\Omega$ is defined so that $k_\Omega(x,y)$ equals the infimum of $k_\Omega$-length($\gamma$) over all rectifiable curves $\gamma \subset \Omega$ joining $x$ with $y$. In case the infimum is attained, the curve for which it happens is called a \emph{quasihyperbolic geodesic}. Any two points of a~proper subdomain of $\R^n$ can be joined by a~quasihyperbolic geodesic, but not necessarily in a~unique manner \cite[Lemma 1]{go79}.

\subsubsection*{Bounded turning}

Let $p\in \mathbb{N}$ and $\Omega\subset\dR$ a domain. We say that $\Omega$ is \emph{homotopically $(p,\varphi)$-bounded turning} (htop $(p, \varphi)$-BT) if every continuous map $f:\mathbb{S}^p\to\Omega$ has a~continuous extension
\begin{equation*}
	g:\closure{\mathbb{B}}^{p+1}\to\Omega \text{ such that } \diam(g(\closure{\mathbb{B}}^{p+1}))\leq \varphi(\diam(f(\mathbb{S}^p))).
\end{equation*}
Observe that for $p = 1$ it boils down to saying that $\Omega$ is simply connected and that any loop in $\Omega \cap B(x,r)$ is contractible in $\Omega \cap B(x, \varphi(r))$ for any $x \in \R^n$ and $r >0$.

In an analoguous manner, we say that $\Omega$ is \emph{homologically $(p, \varphi)$-bounded turning} (hlog $(p, \varphi)$ BT) if each $p$-cycle $z$ in $\Omega$ bounds a $(p+1)$-chain $g$ with $\diam(|g|)\leq \varphi(\diam(|z|))$. Loosely speaking, for $p=1$ it means that each loop forms a~boundary of a~two-dimensional manifold with a~comparable diameter.

In the linear case, the homotopical $p$-dimensional bounded turning condition was firstly introduced by Alestalo and V\"ais\"al\"a in \cite{av94}, whereas its homological counterpart and relationships between them appeared in \cite{av96}. As remarked in the introduction, for $p=0$ both definitions amount to the classical one (any two points $x, y \in \Omega$ can be joined by a~continuum whose diameter does not exceed $C|x-y|$ for some constant $C>1$). It is worth noting that the classical BT property is quantitatively equivalent to LLC-1. 

\subsubsection*{Joinability}
We will define two types of joinability -- homological (\emph{hlog}) for open sets and cohomological (\emph{cohlog}) for closed sets.
Firstly, suppose that $A \subset \dR$ is open. We say $A$ is \emph{hlog outer joinable} if for any $a \in A \setminus \{\infty\}$ and $r>0$, every $p$-cycle from $A \cap B(a,r)$ which bounds a~chain in $A$, bounds a~chain in $A \cap B(a, \varphi(r))$. This property can be equivalently phrased in a~more algebraic (and more concise) manner. Let $A\stackrel{\alpha}{\longrightarrow}B\stackrel{\beta}{\longrightarrow}C$ be a sequence of groups and their homomorphisms. We call it \emph{fast} if $\mathrm{ker}(\beta \alpha) = \mathrm{ker}(\alpha)$. It is evident that hlog outer joinability is equvialent to fastness of the following sequence induced by inclusions
\begin{equation*}
	H_p(A\cap B(a,r))\rightarrow H_p(A\cap B(a,\varphi(r)))\rightarrow H_p(A).
\end{equation*}
On the other hand, we will say that $A$ is \emph{hlog inner joinable} if the following sequence is fast
\begin{equation*}
	H_p(A\backslash \closure{B}(a,\varphi(r)))\rightarrow H_p(A\backslash \closure{B}(a,r))\rightarrow H_p(A).
\end{equation*}

Now assume that $A \subset \dR$ is closed, in which case we will call a~sequence $A\stackrel{\alpha}{\longrightarrow}B\stackrel{\beta}{\longrightarrow}C$ \emph{slow} if $\mathrm{im}(\beta \alpha) = \mathrm{im}(\beta)$. We say that $A$ is \emph{cohlog outer joinable} if the sequence
\begin{equation*}
	H^p(A)\rightarrow H^p(A\cap \closure{B}(a,\varphi(r)))\rightarrow H^p(A\cap\closure{B}(a,r))
\end{equation*}
is slow. Similarly, we say that $A$ is \emph{cohlog inner joinable} if the following sequence is slow
\begin{equation*}
	H^p(A)\rightarrow H^p(A\backslash B(a,r))\rightarrow H^p(A\backslash B(a,\varphi(r))).
\end{equation*}

The above defined joinability is referred to by V\"ais\"al\"a \emph{absolute joinability}, whereas if one checks the corresponding properties for $a \in \R^n$ one gets the definitions of \emph{relative joinability in $\R^n$}. Lemma \ref{lemma: V 2.5} states that if a~set satisfies relative joinability condition with function $\varphi$, it satisfies also this condition for $a \in \R^n$ with function $2\varphi + \mathrm{id}$. In terms of joinability V\"ais\"al\"a formulated Metric Duality Theorem \cite[Theorem 2.7]{v97} whose nonlinear version will be employed in the proof of Theorem \ref{thm main: characterization}. 

\begin{proposition}[Metric Duality Theorem]\label{thm:duality theorem} Suppose that $U$ is an open set in $\dR$ and $p$ is an integer with $0\leq p\leq n-2$. Set $X=\dR\backslash U$ and $q=n-2-p$. Then
	\begin{enumerate}
		\item $U$ is hlog outer $(p,\varphi)$-joinable in $\R^n$ if and only if $X$ is cohlog inner $(q,\varphi)$-joinable in $\R^n$;
		\item $U$ is hlog inner $(p,\varphi)$-joinable in $\R^n$ if and only if $X$ is cohlog outer $(q,\varphi)$-joinable in $\R^n$;
		\item $U$ is hlog $(p,\varphi)$-joinable in $\R^n$ if and only if $X$ is cohlog $(q,\varphi)$-joinable in $\R^n$.
	\end{enumerate}
\end{proposition}
The proof given for $\varphi(r) = Cr$ in \cite[Theorem 2.7]{v97} works for general $\varphi$ without change. Joinability is closely related to bounded turning, as shown in the following

\begin{lemma}\label{lemma: characterization of hlog BT}
	Let $p\in \mathbb{N}$ and $\Omega\subset\dR$ be an open set. Then $\Omega$ is hlog $(p,\varphi)$-BT if and only if $H_p(\Omega)=0$ and $\Omega$ is hlog outer $(p,\varphi)$-joinable, quantitatively.
\end{lemma}
\begin{proof}
	Suppose $\Omega$ is hlog $(p,\varphi)$-BT. Then every $p$-cycle $z$ in $\Omega$ bounds a chain $g\subset\Omega$ such that $\diam(|g|)\leq \varphi(\diam(|z|))$. This implies that for $a\in \Omega \setminus \{\infty\}$ and $r>0$, the sequence
	\begin{equation*} %\label{eq: characterization of hlog BT}
		H_p(\Omega\cap B(a,r))\rightarrow H_p(\Omega\cap B(a, \varphi(2r) + r))\rightarrow H_p(\Omega)
	\end{equation*}
	is fast, which is the definition of $\Omega$ being hlog outer $(p, \varphi(2\cdot) + \mathrm{id})$-joinable. 
	To prove the converse, take a~$p$-cycle $z$ in $\Omega$. It is contained in $B(a, \diam{|z|})$ for some $a \in \Omega \setminus \{\infty\}$. As $\Omega$ is hlog outer $(p, \varphi)$-joinable and $H_p(\Omega) = 0$, the mapping \[H_p(\Omega \cap B(a, \diam{|z|})) \to H_p(\Omega \cap B(a, \varphi(\diam{|z|})))\] is zero as a~homomorphism. Therefore, $z$ bounds a $(p+1)$-chain $g$ such that $\diam{|g|} \leq 2\varphi(\diam{|z|})$, wich means that $\Omega$ is hlog $(p, 2\varphi)$-BT.
\end{proof}

In light of the preceding lemma, it is worth mentioning that htop $(p, \varphi)$-BT property amounts to saying that $\pi_p(\Omega) = 0$ and that a~homotopical counterpart of the hlog outer joinability condition is satisfied, i.e., that the sequence
$$\pi_p(\Omega\cap B(x,r))\rightarrow \pi_p(\Omega\cap B(x,\varphi(r)))\rightarrow \pi_p(\Omega)$$
is fast.

\subsubsection*{Local collars} A domain $\Omega\subset \R^n$ is said to be \emph{locally collared along $\bdary \Omega$} if each point $x\in \bdary \Omega$ has a neighborhood $U$ such that $U\cap \closure{\Omega}$ is homeomorphic to the intersection of a ball and the closure of the upper half space. Any domain whose closure is homeomorphic to a closed ball is locally collared along its boundary. The proof of the lemma below is the same as the proof in the linear case presented in \cite{y09} and thus is omitted here.
\begin{lemma} \cite[Theorems 3.1, 3.2]{y09} \label{lem: local collars}
	Let $\Omega$ be locally collared along $\bdary \Omega$. Then $\Omega$ is $\varphi$-LC-2 if and only if, quantitatively, $\closure{\Omega}$ is $\varphi$-LC-2.
\end{lemma}

\subsubsection*{Uniformity}
A domain $\Omega\subset \R^n$ is \emph{uniform}, if there exists a constant $A\geq 1$ such that each pair of points $x_1,x_2\in\Omega$ can be joined by a rectifiable curve $\gamma$ in $\Omega$ for which
\begin{equation*} \label{eq:definition of uniform}
	\forall\ t \quad \min_{i=1,2} l(\gamma([x_j, \gamma(t)]))\leq Ad(\gamma(t),\bdary\Omega) \quad \text{and} \quad l(\gamma)\leq Ad(x_1,x_2).
\end{equation*}
Such a curve $\gamma$ as above is called an \emph{$A$-uniform curve}. Roughly speaking, uniform curves cannot be twisted too much and have to stay far from the boundary. Bounded uniform domains are John domains and quasiconformal images of uniform domains are typical examples of the more general class of Gromov hyperbolic domains, see~\cite{bhk01}. Recall that this in particular means that uniform domains satisfy the ball separation property \cite{BaloghBuckley}.

\section{Proof of Theorem \ref{thm main: characterization}} \label{sec: main characterization}

\begin{proof}[\textbf{Proof of Theorem~\ref{thm main: characterization}}]
	We will show each of the two equivalencies separately.
	
	$1) \iff 2)$: Proposition \ref{thm:John implies LLC-2} and Lemma \ref{thm: v 3.5} imply that $(0, \varphi)$-John domains are quantitatively $\varphi$-LC-2. Therefore, it suffices to show that if $\Omega$ is $\varphi$-LC-2, then it is $(0,\psi)$-John for some $\psi$ that is quantitatively equivalent to $\varphi$.  We first claim that for each pair of points $x,y\in \Omega$, there exists a quasihyperbolic geodesic $\gamma$ that joins $x$ and $y$ with the property that
	\begin{align}\label{eq:1}
		\text{either }\gamma_{zx} \text{ or } \gamma_{zy} \text{ is contained in } B\Big(z,\psi\big(d(z,\partial \Omega)\big)\Big)
	\end{align}
	for all $z\in \gamma$ with $\psi=\varphi(Ct)$, where $C$     is the constant from the ball separation property. Indeed, if~\eqref{eq:1} fails for some $z\in \gamma$, then there exists $x_0\in \gamma_{zx}$ and $y_0\in \gamma_{zy}$ such that $x_0$ and $y_0$ are outside the ball $B\Big(z,\psi\big(d(z,\partial \Omega)\big)\Big)$. Then, by the \mbox{$\varphi$-LC-2} condition, they can be joined outside $B(z,Cd(z,\partial \Omega))$, which contradicts the fact that $\Omega$ satisfies the ball separation property with constant $C$. Thus~\eqref{eq:1} holds for all $z\in \gamma$, and consequently, $\Omega$ is (0,$\psi$)-John.
	
	$2) \iff 3)$:
	Since $H^1(\closure{\Omega}) = 0$, by the Alexander duality (see e.g.~\cite[Theorem 74.1]{mu84}), $H_{n-2}(U)=0$. Now we can write (bearing in mind that all equivalences are in the quantitative sense) that
	\begin{align*}
		U \text{ is } (n-2,\varphi)&\text{-bounded turning }\\
		&\stackrel{\text{Lemma}~\ref{lemma: characterization of hlog BT}}{\Longleftrightarrow}   U \text{ is hlog outer } (n-2,\varphi)\text{-joinable}\\
		&\stackrel{\text{Proposition}~\ref{thm:duality theorem}}{\Longleftrightarrow}\closure{\Omega} \text{ is cohlog inner } (0,\varphi)\text{-joinable}\\
		&\stackrel{\text{Lemma}~\ref{thm: v 3.10}}{\Longleftrightarrow}\closure{\Omega} \text{ is } \varphi \text{-LC-2}.
	\end{align*}
	Eventually, observe that since $\Omega$ is locally collared along $\partial\Omega$, Lemma \ref{lem: local collars} allows us to conclude that $\Omega$ is $\varphi$-LC-2 if and only if, quantitatively, $\closure{\Omega}$ is $\varphi$-LC-2.
\end{proof}

Observe that at the very beginning of the proof, we infer that $(0,\varphi)$-John domains are $\varphi$-LC-2, quantitatively. It had been previously proved for bounded planar and simply connected domains under some technical assumptions on $\varphi$, see \cite[Theorem 1.3]{gk14}. What is more, Theorem~\ref{thm main: characterization} is sharp in the sense that the ball separation property cannot be dropped as indicated by the following instructive

\begin{example}\label{example:easy one}
	There exists a domain $\Omega\subsetneq \R^3$ locally collared along $\partial\Omega$ with $H^1(\closure{\Omega})=0$ such that $\Omega$ is LLC-2, but not $C$-diam John for any $C\geq 1$. In particular, $\Omega$ fails to have the ball separation property.
\end{example}

\begin{proof}[Construction of Example~\ref{example:easy one}]
	Simply rotate $|y|=(1-x)^2, 0\le x \le 1$ about the $y$-axis to sweep out a cusp domain in the space (see Figure \ref{fig:ex13}).
\end{proof}

	\begin{figure}
		\begin{tikzpicture}[scale=3]
			
			\filldraw[draw=darkgray,fill=lightgray,scale=1,smooth,variable=\x] plot [domain=-1:0] ({\x},{-(1+\x)^2}) plot [domain=0:1]({\x},{-(\x-1)^2}) arc [x radius=1, y radius=0.25, start angle=0,end angle=-540];
			\filldraw[fill=lightgray,scale=1,smooth,variable=\x] plot [domain=-1:0] ({\x},{(1+\x)^2}) plot [domain=0:1]({\x},{(\x-1)^2}) arc [x radius=1, y radius=0.25, start angle=0,end angle=-180];
			\draw[white] (-0.87,0.13) -- (0,1);
			\draw[thick,lightgray,smooth,variable=\x] plot [domain=-1:-0.65] ({\x},{(1+\x)^2});
			\draw[thick,gray,smooth,variable=\x,opacity=0.3] plot [domain=-1:-0.65] ({\x},{(1+\x)^2});
			\draw[thick, dotted,gray,smooth,variable=\x,opacity=0.3] plot [domain=-1:-0.55] ({\x},{-(1+\x)^2});
			\draw[thick,lightgray,smooth,variable=\x] plot [domain=0.65:1] ({\x},{(1-\x)^2});
			%\draw[thick,gray,smooth,variable=\x,opacity=0.3] plot [domain=0.65:1] ({\x},{(1-\x)^2});
			\draw[thick, dotted,gray,smooth,variable=\x,opacity=0.3] plot [domain=0.55:1] ({\x},{-(1-\x)^2});
			
		\end{tikzpicture}
		\caption{Example~\ref{example:easy one}}
		\label{fig:ex13}
		
	\end{figure}

It is worth noting that the proof of Theorem \ref{thm main: characterization} in the linear case under the stronger assumption that $\Omega$ is quasiconformally equivalent to the unit ball (and not merely satisfying the ball separation property) is significantly easier. Indeed, in this case, the equivalence of 1) and 3) can be obtained as a corollary of \cite[Theorem 5.21]{v97}, whereas the equivalence of 1) and 2) can be proved similarly to the planar case as in \cite{nv91}. Moreover, observe that the proof of the equivalence of 1) and 2) in Theorem~\ref{thm main: characterization} implies that quasihyperbolic geodesics are $\varphi$-diam John curves. We formulate this result as a~separate theorem below.

\begin{theorem}\label{prop:quasihyperbolic geodesic are diam John curves}
	Let $\Omega\subset \R^n$ be a bounded domain with the ball separation property. If $\Omega$ is $\varphi$-LC-2, then $\Omega$ is a $\varphi$-diam John domain and quasihyperbolic geodesics are $\varphi$-diam John curves.	
\end{theorem}

Note that Theorem~\ref{prop:quasihyperbolic geodesic are diam John curves} was previously known to hold for bounded Gromov hyperbolic $\varphi$-diam John domains~\cite[Theorem 3.8]{g15}. As a consequence of Theorem~\ref{prop:quasihyperbolic geodesic are diam John curves} and~\cite[Proof of Theorem 3.1 (2)]{g15}, we obtain the following result that generalizes~\cite[Theorem 3.1 (2)]{g15} by removing the a priori stronger Gromov-hyperbolicity assumption.

\begin{theorem}\label{thm:relation of different John}
	Let $\Omega\subset \R^n$ be a bounded domain with the ball separation property. If $\Omega$ is $\varphi$-LC-2 with $\frac{t}{\varphi^{-1}(t)}$ being non-increasing, then $\Omega$ is $\eta$-length John for \begin{equation}\label{eq:diam to length}
		\eta(t)=C\int_0^{\varphi(Ct)}\Big(\frac{s}{\varphi^{-1}(s)}\Big)^{n-1}ds,
	\end{equation}
	provided this integral converges.
	The statement is essentially sharp in the sense that $\eta$ defined in~\eqref{eq:diam to length} is best possible.
\end{theorem}

\section{Proof of Theorem \ref{thm main: example}} \label{sec: main example}

This section is devoted to showing that the interior of a properly constructed Alexander's horned ball satisfies assumptions of Theorem \ref{thm main: example}. In the literature, this example appears usually in the context of its boundary (Alexander's horned sphere); we will follow the description from \cite[Example 2B.2]{h02}. At the end of this section, we also prove quantitative equivalence of hlog and htop bounded turning properties in the case of planar, simply connected domains.

\subsubsection*{Construction of $\Omega$ and basic properties}
The basic building block in the construction of Alexander's horned ball is a~\emph{pair of horns}, which consists of two identical solid tori twisted within each other so that they do not intersect and linked (see Figure \ref{fig:horns}. The construction is inductive and at its every stage rescaled copies of the first pair of horns are used.

\begin{figure}[hbt]
	\begin{tikzpicture}
		\filldraw[white, fill=white] (5.7,0) to (5.7,-0.7) to (5.45,-4.25)to (5.45,-6) to (11.3,-6) to (11,-0.7) to (11,0)--cycle;
		\shadedraw[lightgray!30!white,left color=gray!95!white, right color=lightgray!30!white,opacity=0.9] (5.7,-0.7) arc [start angle=87, end angle=266, x radius=1.3, y radius=1.772] arc [start angle=266, end angle=447, x radius=1.05, y radius=1.772];

		\shadedraw[top color=gray, bottom color=lightgray!30!white] (5.6,-4.244)  arc [start angle=-93, end angle=87, x radius=1.05, y radius=1.772];
		
		\filldraw[thick,left color=gray, right color=lightgray!30!white] (11,-0.7) arc [start angle=93, end angle=-87, x radius=1.3, y radius=1.772] arc [start angle=-87, end angle=-267, x radius=1.05, y radius=1.772];
		
		\shadedraw[gray, thick, left color=lightgray, right color=gray] (11,-4.22) to [out=175,in=-60] (10.4,-3.6) to [out=120, in=50]  (9.7,-3.2) arc [x radius=1.18, y radius=0.85, start angle=-35, end angle=-180] arc [start angle=180, end angle=40, x radius=1.18, y radius=0.75] to [out=-40, in =-145] (11,-0.7) arc [start angle=93, end angle=-87, x radius=1.3, y radius=1.772] to [out=180,in=-5] (11,-4.22);

		\filldraw[gray, thick, fill=white] (8.12, -2.64) arc [x radius=0.85, y radius=0.65, start angle=-140, end angle=-40] arc [x radius=0.83, y radius=0.65, start angle=-40, end angle=-50] arc [x radius=0.83, y radius=0.6, start angle=50, end angle=130];
		
		\shadedraw[thick,gray, left color=gray!80!white, right color=lightgray!20!white,opacity=1] (5.65,-0.68) -- (5.7,-0.7) to [out=-10,in=-128] (6.65,-2.1) arc [start angle=140, end angle=-1, x radius=1.3, y radius=1.15] arc [x radius=0.83, y radius=0.6, start angle=-78, end angle=-114] arc [x radius=0.71, y radius=0.61, start angle=0, end angle=293] arc [x radius=1.18, y radius=0.85, start angle=-128, end angle=-92.5] arc [start angle=-39, end angle=-140, x radius=1.3, y radius=1.15] to [out=128,in=60] (6.4,-3.444) to [out=-120,in=10] (5.58,-4.244) to [out=-170, in=-200] (5.4,-4.22) arc [start angle=265, end angle=84, x radius=1.3, y radius=1.772];

		\shadedraw[thick,gray, top color=gray, bottom color=lightgray!20!white,opacity=1] (5.65,-0.68) -- (5.7,-0.7) to [out=-10,in=-128] (6.65,-2.1) arc [start angle=140, end angle=-1, x radius=1.3, y radius=1.15] arc [x radius=0.83, y radius=0.6, start angle=-78, end angle=-114] arc [x radius=0.71, y radius=0.61, start angle=0, end angle=293] arc [x radius=1.18, y radius=0.85, start angle=-128, end angle=-92.5] arc [start angle=-39, end angle=-140, x radius=1.3, y radius=1.15] to [out=128,in=60] (6.4,-3.444) to [out=-120,in=10] (5.58,-4.244) to [out=-170, in=-200] (5.4,-4.22) arc [start angle=265, end angle=84, x radius=1.3, y radius=1.772];

		%krzywki
		\draw[gray,thick,opacity=0.5] (5.2,-1.3) to [out=-30,in=160] (6.3,-2.4);
		\draw[gray,thick,opacity=0.4] (4.9,-2.8) to [out=15,in=180] (6.25,-2.71);
		\draw[gray,thick,opacity=0.25] (5.3,-3.7) to [out=30,in=190] (6.3,-3.1);
		\draw[darkgray,thick,opacity=0.2] (11.5,-1.8) to [out=200,in=10] (10.3,-2.5);
		\draw[darkgray,thick,opacity=0.2] (11.45,-3.2) to [out=150,in=-10] (10.3,-2.85);
	\end{tikzpicture}
	\centering
	\caption{A pair of `horns'}
	\label{fig:horns}
\end{figure}

The starting point is a~solid torus whose part is taken out (by cutting the torus with two planes), which is simply a~curved solid cylinder. Let us call it $\closure{B}_0$ and note that it obviously is homeomorphic to a~ball. Then, a~pair of horns is glued to both bases of this cylinder so that the whole structure is contained in the initial torus. In the second step, each solid tori in the glued pair of horns is cut so that a~properly rescaled pair of horns can be glued there (again, so as to create a~decreasing family of compact sets). The set before gluing the second generation of pairs of horns is $\closure{B}_1$ and is homeomorphic to $B_0$. At $j$-th stage of the construction, we add $2^{j-1}$ pairs of horns by gluing them to the set $\closure{B}_j$, which has $2^{j-2}$ pairs of horns produced at the preceding stage (including the pairs of generation $j-1$, which are cut).

\begin{figure}[htb]
	\begin{tikzpicture}[scale=0.9]
		\shadedraw[top color=gray, bottom color=lightgray!50!white, thick] (0,0) arc [x radius=8.25, y radius=4.5, start angle=-180, end angle=0] arc [x radius =8.25, y radius=3.5, start angle=0, end angle=180];
		\filldraw[fill=white,thick] (3,0) arc [x radius=7.5, y radius=3, start angle=-135, end angle=-45] arc [x radius=7.5, y radius=3, start angle=-45, end angle=-52] arc [x radius=6.5, y radius=3.5, start angle=45, end angle=135];
		
		\filldraw[white, fill=white] (5.7,0) to (5.7,-0.7) to (5.45,-4.25)to (5.45,-6) to (11.3,-6) to (11,-0.7) to (11,0)--cycle;
		\shadedraw[lightgray!30!white,left color=gray!95!white, right color=lightgray!30!white,opacity=0.9] (5.7,-0.7) arc [start angle=87, end angle=266, x radius=1.3, y radius=1.772] arc [start angle=266, end angle=447, x radius=1.05, y radius=1.772];
		
		\shadedraw[top color=gray, bottom color=lightgray!30!white] (5.6,-4.244)  arc [start angle=-93, end angle=87, x radius=1.05, y radius=1.772];
		
		\filldraw[thick,left color=gray, right color=lightgray!30!white] (11,-0.7) arc [start angle=93, end angle=-87, x radius=1.3, y radius=1.772] arc [start angle=-87, end angle=-267, x radius=1.05, y radius=1.772];
		
		\shadedraw[gray, thick, left color=lightgray, right color=gray] (11,-4.22) to [out=175,in=-60] (10.4,-3.6) to [out=120, in=50]  (9.7,-3.2) arc [x radius=1.18, y radius=0.85, start angle=-35, end angle=-180] arc [start angle=180, end angle=40, x radius=1.18, y radius=0.75] to [out=-40, in =-145] (11,-0.7) arc [start angle=93, end angle=-87, x radius=1.3, y radius=1.772] to [out=180,in=-5] (11,-4.22);

		%uśmiech w prawym rogu
		\filldraw[gray, thick, fill=white] (8.12, -2.64) arc [x radius=0.85, y radius=0.65, start angle=-140, end angle=-40] arc [x radius=0.83, y radius=0.65, start angle=-40, end angle=-50] arc [x radius=0.83, y radius=0.6, start angle=50, end angle=130];

		%lewy róg
		\shadedraw[thick,gray, left color=gray!80!white, right color=lightgray!20!white,opacity=1] (5.65,-0.68) -- (5.7,-0.7) to [out=-10,in=-128] (6.65,-2.1) arc [start angle=140, end angle=-1, x radius=1.3, y radius=1.15] arc [x radius=0.83, y radius=0.6, start angle=-78, end angle=-114] arc [x radius=0.71, y radius=0.61, start angle=0, end angle=293] arc [x radius=1.18, y radius=0.85, start angle=-128, end angle=-92.5] arc [start angle=-39, end angle=-140, x radius=1.3, y radius=1.15] to [out=128,in=60] (6.4,-3.444) to [out=-120,in=10] (5.58,-4.244) to [out=-170, in=-200] (5.4,-4.22) arc [start angle=265, end angle=84, x radius=1.3, y radius=1.772];

		\shadedraw[thick,gray, top color=gray, bottom color=lightgray!20!white,opacity=1] (5.65,-0.68) -- (5.7,-0.7) to [out=-10,in=-128] (6.65,-2.1) arc [start angle=140, end angle=-1, x radius=1.3, y radius=1.15] arc [x radius=0.83, y radius=0.6, start angle=-78, end angle=-114] arc [x radius=0.71, y radius=0.61, start angle=0, end angle=293] arc [x radius=1.18, y radius=0.85, start angle=-128, end angle=-92.5] arc [start angle=-39, end angle=-140, x radius=1.3, y radius=1.15] to [out=128,in=60] (6.4,-3.444) to [out=-120,in=10] (5.58,-4.244) to [out=-170, in=-200] (5.4,-4.22) arc [start angle=265, end angle=84, x radius=1.3, y radius=1.772];
		
		%krzywki
		\draw[gray,thick,opacity=0.5] (5.2,-1.3) to [out=-30,in=160] (6.3,-2.4);
		\draw[gray,thick,opacity=0.4] (4.9,-2.8) to [out=15,in=180] (6.25,-2.71);
		\draw[gray,thick,opacity=0.25] (5.3,-3.7) to [out=30,in=190] (6.3,-3.1);
		\draw[darkgray,thick,opacity=0.2] (11.5,-1.8) to [out=200,in=10] (10.3,-2.5);
		\draw[darkgray,thick,opacity=0.2] (11.45,-3.2) to [out=150,in=-10] (10.3,-2.85);
		
		\begin{scope}[xshift=212,yshift=-75.5,xscale=0.16, yscale=0.2]
			
			% wybielacz
			\filldraw[white, fill=white] (5.7,-0.7) arc [x radius=7.5, y radius=2, start angle=-110.5, end angle=-68.5] arc [start angle=93, end angle=-87, x radius=1.3, y radius=1.772] arc [x radius=8.25, y radius=7, start angle=-70, end angle=-110] arc [start angle=266, end angle=87, x radius=1.3, y radius=1.772] ;
		\end{scope}

		\shade[top color=gray, bottom color=lightgray!20!white] (5.65,-0.68) -- (5.7,-0.7) to [out=-10,in=-128] (6.65,-2.1) arc [start angle=140, end angle=-39, x radius=1.3, y radius=1.15] arc [x radius=1.18, y radius=0.85, start angle=-92.5, end angle=-109]-- (8.43,-2.82)arc [x radius=0.71, y radius=0.61, start angle=0, end angle=293] arc [x radius=1.18, y radius=0.85, start angle=-128, end angle=-92.5] arc [start angle=-39, end angle=-140, x radius=1.3, y radius=1.15] to [out=128,in=60] (6.4,-3.444) to [out=-120,in=10] (5.58,-4.244) to [out=-170, in=-200] (5.4,-4.22) arc [start angle=265, end angle=84, x radius=1.3, y radius=1.772];
		
		\draw[thick,gray] (5.65,-0.68) -- (5.7,-0.7) to [out=-10,in=-128] (6.65,-2.1) arc [start angle=140, end angle=-140, x radius=1.3, y radius=1.15] to [out=128,in=60] (6.4,-3.444) to [out=-120,in=10] (5.58,-4.244) to [out=-170, in=-200] (5.4,-4.22) arc [start angle=265, end angle=84, x radius=1.3, y radius=1.772];
		
		\draw[thick,gray] (8.43,-2.82)arc [x radius=0.71, y radius=0.61, start angle=0, end angle=293];
		
		\begin{scope}[xshift=212,yshift=-75.5,xscale=0.16, yscale=0.2]
			\filldraw[thick,darkgray,left color=gray, right color=lightgray!30!white] (11,-0.7) arc [start angle=93, end angle=-87, x radius=1.3, y radius=1.772] arc [start angle=-87, end angle=-267, x radius=1.05, y radius=1.772];
			
			\shadedraw[gray, thick, left color=lightgray, right color=gray] (11,-4.22) to [out=175,in=-60] (10.4,-3.6) to [out=120, in=50]  (9.7,-3.2) arc [x radius=1.18, y radius=0.85, start angle=-35, end angle=-180] arc [start angle=180, end angle=40, x radius=1.18, y radius=0.75] to [out=-40, in =-145] (11,-0.7) arc [start angle=93, end angle=-87, x radius=1.3, y radius=1.772] to [out=180,in=-5] (11,-4.22);

			%uśmiech w prawym rogu
			\filldraw[gray, thick, fill=white] (8.12, -2.64) arc [x radius=0.85, y radius=0.65, start angle=-140, end angle=-40] arc [x radius=0.83, y radius=0.65, start angle=-40, end angle=-50] arc [x radius=0.83, y radius=0.6, start angle=50, end angle=130];

			%lewy róg
			\shadedraw[thick,gray, left color=gray!80!white, right color=lightgray!20!white,opacity=1] (5.65,-0.68) -- (5.7,-0.7) to [out=-10,in=-128] (6.65,-2.1) arc [start angle=140, end angle=-1, x radius=1.3, y radius=1.15] arc [x radius=0.83, y radius=0.6, start angle=-78, end angle=-114] arc [x radius=0.71, y radius=0.61, start angle=0, end angle=293] arc [x radius=1.18, y radius=0.85, start angle=-128, end angle=-92.5] arc [start angle=-39, end angle=-140, x radius=1.3, y radius=1.15] to [out=128,in=60] (6.4,-3.444) to [out=-120,in=10] (5.58,-4.244) to [out=-170, in=-200] (5.4,-4.22) arc [start angle=265, end angle=84, x radius=1.3, y radius=1.772];

			\shadedraw[thick,gray, top color=gray, bottom color=lightgray!20!white,opacity=1] (5.65,-0.68) -- (5.7,-0.7) to [out=-10,in=-128] (6.65,-2.1) arc [start angle=140, end angle=-1, x radius=1.3, y radius=1.15] arc [x radius=0.83, y radius=0.6, start angle=-78, end angle=-114] arc [x radius=0.71, y radius=0.61, start angle=0, end angle=293] arc [x radius=1.18, y radius=0.85, start angle=-128, end angle=-92.5] arc [start angle=-39, end angle=-140, x radius=1.3, y radius=1.15] to [out=128,in=60] (6.4,-3.444) to [out=-120,in=10] (5.58,-4.244) to [out=-170, in=-200] (5.4,-4.22) arc [start angle=265, end angle=84, x radius=1.3, y radius=1.772];
			
			%krzywki
			\draw[gray,thick,opacity=0.5] (5.2,-1.3) to [out=-30,in=160] (6.3,-2.4);
			\draw[gray,thick,opacity=0.4] (4.9,-2.8) to [out=15,in=180] (6.25,-2.71);
			\draw[gray,thick,opacity=0.25] (5.3,-3.7) to [out=30,in=190] (6.3,-3.1);
			\draw[darkgray,thick,opacity=0.2] (11.5,-1.8) to [out=200,in=10] (10.3,-2.5);
			\draw[darkgray,thick,opacity=0.2] (11.45,-3.2) to [out=150,in=-10] (10.3,-2.85);
			
		\end{scope}
		
		\begin{scope}[rotate=186,xshift=-246,yshift=92,xscale=0.16, yscale=0.16]
			
			% wybielacz
			\fill[fill=white] (5.7,-0.75) arc [x radius=7.5, y radius=7, start angle=-110.5, end angle=-68.5] arc [start angle=93, end angle=-87, x radius=1, y radius=1.9] arc [x radius=8.25, y radius=7, start angle=-70, end angle=-110] arc [start angle=266, end angle=87, x radius=1.2, y radius=1.9] ;
			
			\fill[fill=white] (10,-0.75) circle [radius=0.5];
		\end{scope}
		\filldraw[thick,gray,fill=lightgray!90!gray] (7.79,-2.265) arc [x radius=1.15, y radius=0.85, start angle=140, end angle=120]--(8,-2.3);
		\begin{scope}[rotate=186,xshift=-246,yshift=92,xscale=0.16, yscale=0.16]
			\fill[fill=gray] (11.5,-3.88) circle [radius=0.5];
			
			\filldraw[thick,darkgray,left color=gray, right color=lightgray!30!white] (11,-0.7) arc [start angle=93, end angle=-87, x radius=1.3, y radius=1.772] arc [start angle=-87, end angle=-267, x radius=1.05, y radius=1.772];
			
			\shadedraw[gray, thick, left color=lightgray, right color=gray] (11,-4.22) to [out=175,in=-60] (10.4,-3.6) to [out=120, in=50]  (9.7,-3.2) arc [x radius=1.18, y radius=0.85, start angle=-35, end angle=-180] arc [start angle=180, end angle=40, x radius=1.18, y radius=0.75] to [out=-40, in =-145] (11,-0.7) arc [start angle=93, end angle=-87, x radius=1.3, y radius=1.772] to [out=180,in=-5] (11,-4.22);

			%uśmiech w prawym rogu
			\filldraw[gray, thick, fill=white] (8.12, -2.64) arc [x radius=0.85, y radius=0.65, start angle=-140, end angle=-40] arc [x radius=0.83, y radius=0.65, start angle=-40, end angle=-50] arc [x radius=0.83, y radius=0.6, start angle=50, end angle=130];

			%lewy róg
			\shadedraw[thick,gray, left color=gray!80!white, right color=lightgray!20!white,opacity=1] (5.65,-0.68) -- (5.7,-0.7) to [out=-10,in=-128] (6.65,-2.1) arc [start angle=140, end angle=-1, x radius=1.3, y radius=1.15] arc [x radius=0.83, y radius=0.6, start angle=-78, end angle=-114] arc [x radius=0.71, y radius=0.61, start angle=0, end angle=293] arc [x radius=1.18, y radius=0.85, start angle=-128, end angle=-92.5] arc [start angle=-39, end angle=-140, x radius=1.3, y radius=1.15] to [out=128,in=60] (6.4,-3.444) to [out=-120,in=10] (5.58,-4.244) to [out=-170, in=-200] (5.4,-4.22) arc [start angle=265, end angle=84, x radius=1.3, y radius=1.772];

			\shadedraw[thick,gray, top color=gray, bottom color=lightgray!20!white,opacity=1] (5.65,-0.68) -- (5.7,-0.7) to [out=-10,in=-128] (6.65,-2.1) arc [start angle=140, end angle=-1, x radius=1.3, y radius=1.15] arc [x radius=0.83, y radius=0.6, start angle=-78, end angle=-114] arc [x radius=0.71, y radius=0.61, start angle=0, end angle=293] arc [x radius=1.18, y radius=0.85, start angle=-128, end angle=-92.5] arc [start angle=-39, end angle=-140, x radius=1.3, y radius=1.15] to [out=128,in=60] (6.4,-3.444) to [out=-120,in=10] (5.58,-4.244) to [out=-170, in=-200] (5.4,-4.22) arc [start angle=265, end angle=84, x radius=1.3, y radius=1.772];

			%krzywki
			\draw[gray,thick,opacity=0.5] (5.2,-1.3) to [out=-30,in=160] (6.3,-2.4);
			\draw[gray,thick,opacity=0.4] (4.9,-2.8) to [out=15,in=180] (6.25,-2.71);
			\draw[gray,thick,opacity=0.25] (5.3,-3.7) to [out=30,in=190] (6.3,-3.1);
			\draw[darkgray,thick,opacity=0.2] (11.5,-1.8) to [out=200,in=10] (10.3,-2.5);
			\draw[darkgray,thick,opacity=0.2] (11.45,-3.2) to [out=150,in=-10] (10.3,-2.85);
			
			\fill[fill=white] (6.1,-4.7) circle [radius=0.5];
			
		\end{scope}
		
		\fill[fill=white] (8.8,-3.72) circle [radius=0.2];
		\fill[fill=lightgray!55!white] (8.375,-3.54) circle [radius=0.03];
		\draw[gray,thick,opacity=0.5] (5.2,-1.3) to [out=-30,in=160] (6.3,-2.4);
		\draw[gray,thick,opacity=0.4] (4.9,-2.8) to [out=15,in=180] (6.25,-2.71);
		\draw[gray,thick,opacity=0.25] (5.3,-3.7) to [out=30,in=190] (6.3,-3.1);
		\fill[fill=white] (7.35,-1.45) circle [radius=0.4];
	\end{tikzpicture}
	\caption{Alexander's horned ball (2nd stage of the construction)}
	\label{fig:alexander}
\end{figure}

\begin{lemma} \label{lem: homeomorphic to ball}
	The set $\closure{\Omega}$ is homeomorphic to the unit ball.
\end{lemma}

\begin{proof}
	There exist a sequence of homeomorphisms $h_j: \closure{B}_{j-1} \to \closure{B}_j$, where $h_0: \closure{\mathbb{B}} \to \closure{B}_0$. Creating $\closure{B}_j$ from $\closure{B}_{j-1}$ affects only a neighbourhood of $\closure{B}_j \setminus \closure{B}_{j-1}$ and hence the sequence $f_j = h_j \circ \dots \circ h_0$ converges uniformly to a~continuous map $f$. Since $f$ is also injective, $f$ is a~homoemorphism, as $\closure{\mathbb{B}}$ and $f(\closure{\mathbb{B}})$ are compact.
\end{proof}

\begin{lemma} 
	The set $U:=\R^3 \setminus \closure{\Omega}$ has perfect fundamental group.
\end{lemma}
The proof of this fact can be found in \cite[p.171]{h02}
% To prove that it is non-trivial, one could only say that any loop around a torus in the construction of AHB cannot be contracted.

\subsubsection*{Uniformity}

For simplicity, we now will use the word 'torus' to indicate a~solid torus. It is convenient to think of this construction in terms of a~binary tree, whose root is the first cut torus denoted with $T_0$ (previously it was called $\closure{B}_0$) and whose two children are the tori in the first pair of horns. Then, each torus has again two children. All tori, except for the starting one, are self-similar. We will say that a~point $x$ belongs to a~torus constructed at stage $j = 1, 2, \ldots$ and write $x \in T_j$ if it belongs to one of the $4^{j-1}$ tori produced at stage $j$ and does not belong to any torus produced at a~later stage. This naturally means that such a~point belongs to some torus with a~cylinder cut out.

Let us also introduce some notation used in the proof of uniformity. Let $d_j$ denote the Euclidean distance between two horns within one pair constructed at $j$-th stage; $r_j$ radius of the small circle of torus $T_j$ and $p_j$ the length of the center line of a~torus $T_j$. By a~\emph{center line} of a~(solid) torus we mean the circle inside the torus, whose distance to the boundary at every point equals the radius of the smaller circle. Given torus $T_j$, we construct $T_{j+1}$ simply by rescaling $T_j$ by a~factor $\lambda$. Therefore, $d_j = \lambda^{j-1}d_1$, $r_j = \lambda^{j-1}r_1$ and $p_j = \lambda^{j-1}p_1$.

\begin{proposition} \label{prop: uniformity}
	The domain $\Omega$ is uniform.
\end{proposition}

\begin{proof}
	Take any two points $x, y \in \Omega$, we will prove that there exists a~$D$-uniform curve $\gamma$ which connects these points. Firstly, assume that $x, y \in T_0$. Since a~curved cylinder $T_0$ is a~bounded Lispchitz domain, it is also uniform with some constant, say $D_1$. Therefore, a~$D_1$-uniform curve connecting $x$ and $y$ clearly exists. Similarly, if both points lie in the same torus $T_j$ for $j \geq 1$, one should observe that a~torus with a cylinder cut out is also uniform (as it also is a~bounded Lipschitz domain). Since all such tori are just rescaled copies of the first ones, they are uniform with the same constant, say $D_2$.
	
	In the case when the chosen points belong to different tori, we can assume without loss of generality that $x \in T_k$ whereas $y \in T_\ell$ and $k \geq \ell$. Observe that in order to connect points from different tori with a~curve, one has to find the torus which is their lowest common ancestor in our binary tree. Let it be a~torus $T_n$ for some $n \in \{0, 1, \ldots, \ell-1\}$. The curve will then wind up from $x$ to $T_n$ and then descend to $y$. In order to construct a uniform curve, we will use center lines of the tori through which the curve goes.
	
	Given $x$ and $y$, let $x'$ and $y'$ denote the projections of $x$ and $y$, respectively, onto the suitable center lines. Choose a~point $z \in T_n$ lying on the center line of $T_n$. Connect $x$ and $x'$ with a~segment, do the same with $y$ and $y'$ and then use center lines of the appropriate tori to connect $x'$ to $z$ and $y'$ to $z$. We will now prove that thus constructed curve is indeed uniform and that its constant does not depend on $k, \ell$ and $n$.
	
	The length of $\gamma(x, z) \cap T_j$ for $j = k+1, \ldots, n$ cannot exceed the length of the whole center line, i.e $p_j$. To accommodate for the line segment between $x$ and $x'$, we can write that $l \left(\gamma \cap T_k\right) \leq 2p_k$, which yields the following estimate
	\begin{align*}
		l(\gamma(x,z)) &\leq p_n + p_{n+1} + \ldots + 2p_k \leq 2(p_n + p_{n+1} + \ldots +p_k) = 2p_1(\lambda^{n-1} + \ldots + \lambda^{k-1}).
	\end{align*}
	Naturally, this quantity bounds the length of $\gamma(z, y)$ as well, which means that
	\begin{equation*}
		l(\gamma) \leq 4\frac{\lambda^{n-1}}{1 - \lambda}p_1.
	\end{equation*}
	It is eminent from the construction that the Euclidean distance between $x$ and $y$, whose lowest common ancestor belongs to $T_n$ must be at least $d_{n+1} = \lambda^{n}d_1$. Therefore, $l(\gamma) \leq D d(x, y)$ for any $D \geq \frac{4p_1}{\lambda(1 - \lambda)d_1}$. It is worth noting that if the lowest common ancestor of $x$ and $y$ is $T_0$, one gets that $D$ needs to satisfy
	\begin{equation*}
		4 \left( p_0 + \frac{\lambda^{n-1}}{1 - \lambda}p_1 \right) \leq Dd_1.
	\end{equation*}
	
	Let us now turn to the second condition in the definition of uniformity. To begin with, fix $t$ such that $\gamma(t) \in \gamma(x, z)$ (for $t$ such that $\gamma(t) \in \gamma(z, y)$ an analogous argument works). If $\gamma(t)$ lies on the line segment connecting $x$ and $x'$, the second condition trivially holds with $D = 1$, since the distance to the boundary is increasing along this segment. If $\gamma(t)$ does not lie on this segment, it has to lie on a~center line inside a~torus $T_j$ for some $j \in \{n, n+1, \ldots, k\}$. Then, $d(\gamma(t), \partial \Omega) = r_j$, whereas 
	\begin{equation*}
		l(\gamma(x, \gamma(t))) \leq 2 \frac{\lambda^{j-1}}{1 - \lambda} p_1.
	\end{equation*}
	Consequently, the condition is satisfied with any $D \geq \frac{2p_1}{(1 - \lambda)r_1}$. Just like before, a slight amendment is necessary in case if $n$ equals zero.
	
	All in all, one can find a~constant $D_3$ such that any two points belonging to different tori of the interior of the Alexander's horned sphere can be joined by a~$D_3$-uniform curve. Choosing $D = \max{ \{D_1, D_2, D_3\} }$ concludes the proof.
\end{proof}

\begin{proof}[\textbf{Proof of Theorem \ref{thm main: example}}]
	Lemma \ref{lem: homeomorphic to ball} and Proposition \ref{prop: uniformity} show that the first two conditions are satisfied. The fundamental group of $U$ is non-trivial, hence $U$ cannot be homotopically $(1,C)$-bounded turning for any $C$. However, application of Theorem \ref{thm main: characterization} implies that $U$ satisfies the homological bounded turning condition.
\end{proof}

It is worth noting that to obtain a domain in $\R^3$ which is hlog $(1,C)$-bounded turning and not htop $(1,C)$-bounded turning, the infinite construction is in fact needed as the complement of each finite stage of this construction is not homologically trivial. Complement of Alexander's horned ball provides an example that for $n=3$ homological and homotopical bounded turning properties do not coincide, which is in turn true for $n=2$, as shown in Proposition \ref{prop: planar bt}. This result, though hardly surprising and probably known to specialists, is of separate interest  and we have decided to include the proof for completeness of the presentation. Note also that this provides a partial answer to \cite[Open problem 4.8]{av94} in the simplest case $n=2$.

\begin{lemma} \label{lemma: Jordan boundary}
	Let $\gamma$ be a closed curve (not necessarily Jordan) contained in a~simply connected domain $\Omega \subset \R^2$. For any $\eps >0$ there exists a~compact set $K \subset \Omega$ containing $\gamma$ whose boundary is a~Jordan curve and for which it is true that
	\[
	\diam{K} \leq (1 + \eps) \diam{\gamma}.
	\]
\end{lemma}

\begin{proof}
	Firstly, let us observe that there exists exactly one unbounded connected component of $\R^2 \setminus \gamma$. Indeed, if $\Gamma: \Sf^1 \to \Omega$ is a~parametrization of $\gamma$, then the function $|\Gamma|: \Sf^1 \to \R$ is s a continuous function over a compact set and hence attains a~maximal value $R$. Therefore, all points outside the disk $B(0,R)$ belong to an unbounded connected component, which has to be unique since connected components are disjoint.
	
	Choose $\delta>0$ so that
	\begin{equation*}
		\sqrt{2} \delta < \min{ \left \{ \mathrm{dist}(\gamma, \partial \Omega), \, \eps/4 \diam{\gamma} \right \}}
	\end{equation*}
	and consider the family of closed squares $Q(x_i, \delta/2)$ of sidelengths $\delta$ whose centers $x_i$ belong to $(\delta \mathbb{Z})^2$. There are those who intersect $\gamma$ (\emph{intersection} squares) and those who do not. Within the latter group there are \emph{outer} ones -- points from these squares can be connected with some point of $S(0, R)$ without crossing an intersection square and \emph{inner} ones -- without this property. Let us choose all the intersection and inner squares (since $\gamma$ is contained in the ball $B(0, R)$, there have to be finitely many of them) to get a~simply connected set $K$. Indeed, were it not simply connected, there would be an outer cube surrounded by intersection or inner cubes, which is impossible. The boundary of $K$ consists of sides of squares, which clearly forms a~Jordan curve. 
	
	What is more, a~side of an inner square cannot be contained in $\partial K$ as then it would be possible to connect points from such a~square with $S(0,R)$, which clearly lies in the unbounded component of $\R^2 \setminus \gamma$. Therefore, for any points $x, y \in \partial K$ there exist points $z_1, z_2 \in \gamma$ such that $|x - z_1| < \delta \sqrt{2}$ and $|y - z_2| < \delta \sqrt{2}$. Triangle inequality and the choice of $\delta$ allows to conclude that
	\begin{equation*}
		\diam{K} = \diam{\partial K} \leq 2\sqrt{2} \delta + \diam{\gamma} \leq (1 + \eps)\diam{\gamma}.
	\end{equation*}
	Clearly, the choice of $\delta$ ensures also that $K \subset \Omega$.
\end{proof}

\begin{lemma} \label{lemma: planar htop bt}
	Let $\Omega \subset \R^2$ be a simply connected domain. Then $\Omega$ is $(1, 1+\eps)$-bounded turning for any $\eps > 0$.
\end{lemma}

\begin{proof}
	Throughout the proof $\closure{\B}^2$ will be denoted as $\D$. We need to show that any continuous function $\Gamma: \mathbb{S}^1 \to \Omega$ can be extended onto the disk to a continuous $\tilde{\Gamma}: \D \to \Omega$ such that \[\diam{(\tilde{\Gamma}(\D))} \leq (1+ \eps)\diam{(\Gamma(\Sf^1))}.\]
	
	Let $\Gamma(\Sf^1) := \gamma$ and choose set $K$ from Lemma \ref{lemma: Jordan boundary} and the homeomorphism $G: \D \xrightarrow{\text{onto}} K$ granted by the Jordan-Schoenflies theorem. The mapping $G^{-1} \circ \Gamma: \Sf^1 \to \D$ is well defined (as $\gamma \subset K$) and can be easily extended onto the disk by
	\[
	\widetilde{G^{-1} \circ \Gamma}(x) = |x| (G^{-1} \circ \Gamma)(x/|x|).
	\]
	Take $\tilde{\Gamma} := G \circ \widetilde{G^{-1} \circ \Gamma}: \D \to K$ to be the requested extension. Indeed, it is continuous as a~composition of continuous mappings and for any $x \in \Sf^1$ we have $\tilde{\Gamma}(x) = \Gamma(x)$. Eventually, due to the fact that $\tilde{\Gamma}(\Sf^1) = K$ and the properties of $K$, one gets that
	\[
	\diam{\tilde{\Gamma}(\Sf^1)} \leq (1+\eps) \diam{\gamma} = (1+ \eps) \diam{\Gamma(\Sf^1)}.
	\]
	Geometrically, the idea is to find a~Jordan domain containing the curve $\gamma$ and contraction of its boundary immediately gives contraction of $\gamma$ as well thanks to Jordan-Schoenflies theorem. The extension $\tilde{\Gamma}$ can be therefore thought of as the homotopy between $\gamma$ and a~point.
\end{proof}

\begin{proposition} \label{prop: planar bt}
	Let $\Omega \subset \R^2$ be a domain. Then $\Omega$ is htop $(1,C)$-bounded turning if and only if it is hlog $(1, C')$-bounded turning.
\end{proposition}

\begin{proof}
	Clearly, if $\Omega$ is homotopically $(1,C)$-bounded turning, then it also satisfies the homological counterpart of this condition. The other direction needs invoking a~few algebraic topology facts. Firstly, it is known that $\pi_1(\Omega)$ is a~free group (\cite[Section 4.2.2]{Stillwell}). Hence its abelianization, $H_1(\Omega)$, is a free abelian group, i.e., a~direct product of the same number of $\mathbb{Z}$ as the number of generators of $\pi_1(\Omega)$. Therefore, if $\pi_1(\Omega)$ were non-trivial, so would $H_1(\Omega)$. Assume now $\Omega$ to be hlog $(1,C)$-bounded turning. It implies that $H_1(\Omega) = 0$ and, in view of the previous argument, it also means that $\pi_1(\Omega) = 0$. Therefore, by Lemma \ref{lemma: planar htop bt}, we can conclude that $\Omega$ is also htop $(1,C)$-bounded turning for some $C > 1$.
\end{proof}

\setcounter{theorem}{0}
\renewcommand{\thetheorem}{A.\arabic{theorem}}

\setcounter{theorem}{0}
\renewcommand{\thetheorem}{A.\arabic{theorem}}
\section*{Appendix}\label{sec: appendix}

In this appendix, nonlinear versions of joinability properties are included for convenience of the reader since they are used in the proof of Theorem \ref{thm main: characterization}. The proofs in the generalized case are essentially the same as in the original paper by V\"ais\"al\"a~\cite{v97} and can be obtained by replacing $cr$ with $\varphi(r)$. We include proofs of Lemma \ref{lemma: V 2.5} and Proposition \ref{thm:John implies LLC-2} in order to show the original reasoning and how quickly it can be adapted to capture the nonlinear setting.

\begin{lemma}\cite[Lemma 2.5]{v97}\label{lemma: V 2.5}
	If $p\geq 0$ and if $A\subset \dR$ is hlog outer $(p,\varphi)$-joinable for all $a\in A\backslash \{\infty\}$ and $r>0$, then $A$ is hlog outer $(p,2\varphi+\mathrm{id})$-joinable in $\R^n$. The same is true for the other three joinability properties as well.
\end{lemma}
\begin{proof} 
	Let $a\in\R^n$, $r>0$. Writing $\varphi'(r)=2\varphi(r)+r$, we must show that the sequence
	\begin{equation*}
		H_p(A\cap B(a,r))\rightarrow H_p(A\cap B(a,\varphi'(r)))\rightarrow H_p(A)
	\end{equation*}
	is fast. If $A\cap B(a,r)=\emptyset$, the first group is trivial, and, consequently, the sequence is fast. If $A\cap B(a,r)\neq \emptyset$, choose a point $x\in A\cap B(a,r)$. Now
	\begin{equation*}
		B(a,r)\subset B(x,2r)\subset B(x,2\varphi(r))\subset B(a,\varphi'(r)),
	\end{equation*}
	and we obtain the commutative diagram
	\[
	\begin{CD}
		H_p(A\cap B(a,r)) @>>> H_p(A\cap B(a,\varphi'(r)))@>>>H_p(A) \\
		@VVV @AAA  \\
		H_p(A\cap B(x,2r)) @>>> H_p(A\cap B(x,2\varphi(r))) @>>>H_p(A).
	\end{CD}
	\]
	Since the lower row is fast, so is the upper row.
	
	Next assume that $A$ is hlog inner $(p,\varphi)$-joinable and that $a\in M$, $r>0$. We must show that the sequence
	\begin{equation*}
		H_p(A\backslash \closure{B}(a,\varphi'(r)))\rightarrow H_p(A\backslash \closure{B}(a,r))\rightarrow H_p(A)
	\end{equation*}
	is fast. If $A\cap \closure{B}(a,r)=\emptyset$, the sequence is trivially fast since then the second map is the identity. If $A\cap \closure{B}(a,r)\neq\emptyset$, choose $x\in A\cap \closure{B}(a,r)$- Now
	\begin{equation*}
		\closure{B}(a,r)\subset \closure{B}(x,2r)\subset \closure{B}(x,2\varphi(r))\subset \closure{B}(a,\varphi'(r)),
	\end{equation*}
	and we can proceed essentially as in the first case.
	
	The cohlog cases are treated by analogous arguments.
\end{proof}

\begin{lemma}\cite[Theorem 3.5]{v97}\label{thm: v 3.5}
	Let $A\subset\dR$. Then $A$ is hlog outer $(0,\varphi)$-joinable in $\R^n$ if and only if every path component of $A$ is pathwise $\varphi$-LC-1. Moreover, $A$ is hlog inner $(0,\varphi)$-joinable in $\R^n$ if and only if every path component of $A$ is pathwise $\varphi$-LC-2.
\end{lemma}

\begin{lemma}\cite[Theorem 3.10]{v97}\label{thm: v 3.10}
	Let $A\subset \dR$ be compact. Then $A$ is cohlog outer $(0,\varphi)$-joinable in $\R^n$ if and only if every path component of $A$ is continuumwise $\varphi$-LC-1, quantitatively. Moreover, $A$ is cohlog inner $(0,\varphi)$-joinable in $\R^n$ if and only if every path component of $A$ is continuumwise $\varphi$-LC-2, quantitatively.
\end{lemma}

\begin{lemma}\cite[Lemma 5.3]{v97}\label{lemma:isomorphism}
	Suppose that $0\leq p\leq n-2$ and that $A,V\subset\dR$ are such that $V$ is open, $H_p(V)=0$, and $V^c\subset int A$. Then the map
	$H_p(A\cap V)\to H_p(A)$ is an isomorphism.
\end{lemma}

\begin{proposition} \cite[Lemma 5.5]{v97} \label{thm:John implies LLC-2}
	Let $U\subset\dR$ be a $(p,\varphi)$-John domain with $0\leq p\leq n-2$. Then $U$ is hlog inner $(p,\varphi)$-joinable, quantitatively.
\end{proposition}
\begin{proof}
	Let $a \in \R^n$ and $r>0$. Let $z$ be a $p$-cycle in $U\backslash\closure{B}(a,2\varphi(r)+r)$ bounding in $U$. We need to show that $z$ bounds in $U\backslash\closure{B}(a,r)$. Since $U$ is $(p,\varphi)$-John, $z=\bdary g$ for some $(p+1)$-chain $g$ satisfying condition \eqref{eq: general John}. If $|g|\cap \closure{B}(a,r)=\emptyset$, there is nothing to prove. In the opposite case, we fix a point $x\in |g|\cap \closure{B}(a,r)$. Now
	\begin{equation*}
		2\varphi(r)=2\varphi(r)+r-r< d(x,|z|)\leq \varphi(d(x,\bdary U)),
	\end{equation*}
	and hence $\closure{B}(a,r)\subset \closure{B}(x,\varphi^{-1}(2\varphi(r)))\subset U$. Applying Lemma~\ref{lemma:isomorphism} with $A=U$ and $V=\closure{B}(a,r)^c$ we see that $z$ bounds in $U\backslash\closure{B}(a,r)$.
\end{proof}

\end{document}